\documentclass[a4paper,11pt]{amsart}

\usepackage{graphicx}
\usepackage{mathptmx}
\usepackage{amsmath}
\usepackage{amssymb}
\usepackage{enumitem}
\usepackage{xcolor}
\usepackage{pgfplots}

\newmuskip\pFqmuskip

\newcommand*\pFq[6][8]{%
  \begingroup 
  \pFqmuskip=#1mu\relax
  \mathcode`=\string"8000
  \begingroup\lccode`\~=`\,
  \lowercase{\endgroup\let~}\pFqcomma
  F^{#2}_{#3}{\left(\genfrac..{0pt}{}{#4}{#5}\bigg|#6\right)}%
  \endgroup
}
\newcommand{\pFqcomma}{\mskip\pFqmuskip}

\newtheorem{theorem}{Theorem}[section]

\newtheorem{remark}[theorem]{Remark}

\begin{document}

\title[ ]{Degenerate Euler-Seidel Method for Degenerate Bernoulli, Euler, and Genocchi Polynomials}

\author{Taekyun  Kim}
\address{Department of Mathematics, Kwangwoon University, Seoul 139-701, Republic of Korea}
\email{tkkim@kw.ac.kr}

\author{Dae San  Kim}
\address{Department of Mathematics, Sogang University, Seoul 121-742, Republic of Korea}
\email{dskim@sogang.ac.kr}

\author{Hyunseok Lee}
\address{Department of Mathematics, Kwangwoon University, Seoul 139-701, Republic of Korea}
\email{luciasconstant@gmail.com}

\author{Kyo-Shin Hwang}
\address{Graduate School of Education, Yeungnam University, Gyeongsan 3854, Republic of Korea}
\email{kshwang@yu.ac.kr}


\begin{abstract}
This paper introduces a degenerate version of the Euler-Seidel method by incorporating a parameter $\lambda$ into the classical recurrence relation. We define a degenerate Euler-Seidel matrix associated with an initial sequence and establish corresponding $\lambda$-generalized binomial identities and generating function relations. By applying this method to the degenerate Bernoulli, Euler, and Genocchi polynomials, we derive several new combinatorial identities. This work extends the classical Euler-Seidel method to the domain of degenerate special polynomials and numbers, providing a new framework for studying their properties.
\end{abstract}

\maketitle

\markboth{\centerline{\scriptsize Degenerate Euler-Seidel Method for Degenerate Bernoulli, Euler, and Genocchi Polynomials}}
{\centerline{\scriptsize Taekyun Kim, Dae San Kim, Hyunseok Lee and Kyo-Shin Hwang}}

\section{Introduction}

For any nonzero $\lambda \in \mathbb{R}$, the degenerate exponentials are defined by
\begin{equation}\label{1}
e_\lambda^x(t) = \sum_{n=0}^{\infty} (x)_{n,\lambda}\frac{t^n}{n!}, \qquad e_\lambda(t)=e_\lambda^{1}(t),
\quad (\text{see } [13,15-19]),
\end{equation}
where the degenerate falling factorial sequence is given by
\begin{equation}
(x)_{0,\lambda} = 1, \quad
(x)_{n,\lambda} = x(x-\lambda)(x-2\lambda)\cdots(x-(n-1)\lambda),
\quad (n \ge 1). \label{2}
\end{equation}
In [2,3], Carlitz introduced the degenerate Bernoulli polynomials given by
\begin{equation}\label{3}
\frac{t}{e_\lambda(t)-1}\, e_\lambda^x(t)
= \sum_{n=0}^{\infty} \beta_{n,\lambda}(x)\frac{t^n}{n!},\ \ (\textrm{see}\ [5,13,18,24]).
\end{equation}
When $x=0$,\ \ $\beta_{n,\lambda}= \beta_{n,\lambda}(0),\ (n \ge 0)$, are called the degenerate Bernoulli numbers.
From \eqref{3}, we note that
\begin{equation}\label{4}
\beta_{n,\lambda}(x) = \sum_{k=0}^{n} \binom{n}{k}\beta_{k,\lambda} (x)_{n-k,\lambda}, \ (n \ge 0).
\end{equation}
In addition, the degenerate Euler polynomials are defined by
\begin{equation}\label{5}
\frac{2}{e_\lambda(t)+1}\, e_\lambda^x(t)
= \sum_{n=0}^{\infty} \mathcal{E}_{n,\lambda}(x)\frac{t^n}{n!},\ \ (\textrm{see}\ [5,8,13,18,22]).
\end{equation}
When $x=0$,\ \  $\mathcal{E}_{n,\lambda} = \mathcal{E}_{n,\lambda}(0), \  (n \ge 0)$, are called the  degenerate Euler numbers.
By \eqref{5}, we get
\begin{equation}\label{6}
\mathcal{E}_{n,\lambda}(x) = \sum_{k=0}^{n} \binom{n}{k} \mathcal{E}_{k,\lambda}(x)_{n-k,\lambda}.
\end{equation}
The degenerate Genocchi polynomials are defined by
\begin{equation}\label{7}
\frac{2t}{e_\lambda(t)+1}\, e_\lambda^x(t) = \sum_{n=0}^{\infty} \mathcal{G}_{n,\lambda}(x)\frac{t^n}{n!},\ \ (\textrm{see}\ [8,19,22]).
\end{equation}
When $x=0$,\ \  $\mathcal{G}_{n,\lambda} = \mathcal{G}_{n,\lambda}(0), \ (n \ge 0)$, are called the  degenerate Genocchi numbers.
Note that
\begin{equation}\label{8}
\mathcal{G}_{0,\lambda} = 0, \quad \mathcal{G}_{n,\lambda}(x)= \sum_{k=0}^{n} \binom{n}{k} \mathcal{G}_{k,\lambda}\,
(x)_{n-k,\lambda},\quad (n \ge 1).
\end{equation}
From \eqref{7}, we have
\begin{equation}\label{9}
\begin{split}
\sum_{n=0}^{\infty} \mathcal{G}_{n,\lambda}\frac{t^n}{n!}
&= \frac{2t}{e_\lambda(t)+1} = 2\left(\frac{t}{e_\lambda(t)-1}- \frac{2t}{e_{\lambda}^{2}(t)-1}\right) \\
&= 2\left(\frac{t}{e_\lambda(t)-1}- \frac{2t}{e_{\lambda/2}(2t)-1}\right)
= 2 \sum_{n=0}^{\infty}\left(\beta_{n,\lambda}- 2^{n}\beta_{n,\lambda/2}\right)\frac{t^n}{n!}.
\end{split}
\end{equation}
Thus, by \eqref{9}, we obtain
\begin{equation}\label{10}
\mathcal{G}_{n,\lambda}= 2\left(\beta_{n,\lambda}- 2^{n}\beta_{n,\lambda/2}\right),\quad (n \ge 0).
\end{equation}

For a given sequence $(a_{n})_{n \ge 0}$, the Euler-Seidel matrix associated with this sequence is determined recursively by
\begin{equation}\label{11}
\begin{aligned}
a_{0,n} &= a_n, \quad (n \ge 0),\\
a_{k,n} & = a_{k-1,n} + a_{k-1,n+1}, \quad (n \ge 0,\; k \ge 1),\quad (\textrm{see}\ [1,7,9,13,20]).
\end{aligned}
\end{equation}
The Euler-Seidel matrix $(a_{n,k})_{n,k \ge 0}$ associated with the sequence $(a_{n})_{n \ge 0}$ is given by
\begin{equation}\label{12}
(a_{n,k})_{n,k \ge 0} =
\begin{pmatrix}
a_{0,0} & a_{0,1} & a_{0,2} & a_{0,3} & \cdots \\
a_{1,0} & a_{1,1} & a_{1,2} & a_{1,3} & \cdots \\
a_{2,0} & a_{2,1} & a_{2,2} & a_{2,3} & \cdots \\
\vdots  & \vdots  & \vdots  & \vdots
\end{pmatrix}.
\end{equation}
From \eqref{11}, we get the binomial identities:
\begin{equation}\label{13}
 a_{n,0}= \sum_{k=0}^{n} \binom{n}{k} a_{0,k}\quad\textrm{and}\quad\
a_{0,n} = \sum_{k=0}^{n} \binom{n}{k} (-1)^{\,n-k} a_{k,0},
\quad(\textrm{see}\ [9]).
\end{equation}
Let $A(t)=\sum_{n=0}^{\infty}a_{0,n}\frac{t^{n}}{n!}$ be the generating function of the initial sequence $(a_{0,n})_{n\ge 0}$. Then the generating function $\overline{A}(t)$ of the final sequence $(a_{n,0})_{n\ge 0}$ is given by
\begin{equation}
\overline{A}(t)=\sum_{n=0}^{\infty}a_{n,0}\frac{t^{n}}{n!}=e^{t}A(t),\quad (\mathrm{see}\ [9]).\label{14}
\end{equation} \par
The purpose of this paper is to extend these results by studying a degenerate version of the Euler-Seidel method. This generalization introduces a parameter $\lambda$ into the recurrence relation. The results are applied to study and derive new combinatorial identities for sequences like the degenerate Bernoulli, Euler, and Genocchi polynomials.

For a given sequence $\big(a_{n}(x|\lambda)\big)_{n \ge 0}$, we consider the degenerate Euler-Seidel matrix associated with this sequence which is determined recursively by (see \eqref{11})
\begin{equation*}
\begin{aligned}
&a_{0,n}(x\mid\lambda)=a_n(x\mid\lambda), \quad (n \ge 0), \\
&a_{k,n}(x\mid\lambda) = \bigl(1-(k-n)\lambda\bigr)\, a_{k-1,n}(x\mid\lambda)
+ a_{k-1,n+1}(x\mid\lambda), \quad (n \ge 0,\, k \ge 1).
\end{aligned}
\end{equation*}
The degenerate Euler--Seidel matrix associated with  $\big(a_{n}(x\mid\lambda)\big)_{n\ge 0}$
  is given by
\begin{equation}\label{15}
\bigl(a_{n,k}{(x\mid\lambda)}\bigr)_{n,k \ge 0}=
\begin{pmatrix}
a_{0,0}{(x\mid\lambda)} & a_{0,1}{(x\mid\lambda)} & a_{0,2}{(x\mid\lambda)} & \cdots \\
a_{1,0}{(x\mid\lambda)} & a_{1,1}{(x\mid\lambda)} & a_{1,2}{(x\mid\lambda)} & \cdots \\
a_{2,0}{(x\mid\lambda)} & a_{2,1}{(x\mid\lambda)} & a_{2,2}{(x\mid\lambda)} & \cdots \\
\vdots & \vdots & \vdots
\end{pmatrix}.
\end{equation}

The sequence $\big(a_{0,n}(x|\lambda)\big)_{n \ge 0}$ is the {\it{initial degenerate sequence}}, and $\big(a_{n,0}(x|\lambda)\big)_{n\ge 0}$ is the {\it{final degenerate sequence}}. The following $\lambda$-generalized binomial identities are established (see Theorems 3.1) using the generalized falling and rising factorials, $(1-\lambda)_{n-k,\lambda}$ and $\langle1-\lambda \rangle_{n-k,\lambda}$, (see \eqref{2}) for the degenerate case:
\begin{align*}
&a_{n,0}{(x\mid\lambda)}=\sum_{k=0}^{n}\binom{n}{k}(1-\lambda)_{\,n-k,\lambda}\,a_{0,k}{(x\mid\lambda)}, \\
&a_{0,n}{(x\mid\lambda)}=\sum_{k=0}^{n}\binom{n}{k}(-1)^{\,n-k}\langle 1-\lambda\rangle_{\,n-k,\lambda}\,
a_{k,0}{(x\mid\lambda)},
\end{align*}
where the degenerate rising factorial sequence is given by
$$
\langle x\rangle_{0,\lambda}=1, \ \ \langle x\rangle_{n,\lambda}= x(x+\lambda)(x+2\lambda)\cdots(x+(n-1)\lambda),
\quad (n\ge1).$$
These correspond to the classical binomial identities in \eqref{13}, and they lead to a degenerate version of Seidel's formula for the generating functions
\begin{equation*}
\overline{A}_\lambda(x,t)=e_\lambda^{\,1-\lambda}(t)\,A_\lambda(x,t),
\end{equation*}
where $A_{\lambda}(x,t)=\sum_{n=0}^{\infty}a_{0,n}(x|\lambda)\frac{t^{n}}{n!}$ and $\overline{A}_\lambda(x,t)=\sum_{n=0}^{\infty} a_{n,0}{(x\mid\lambda)}\frac{t^n}{n!}$ (see Theorem 3.2).
The results derived from the degenerate Euler-Seidel method are applied to the degenerate Bernoulli, Euler, and Genocchi polynomial sequences, yielding various combinatorial identities (see Theorems 3.3-3.5): for $n \ge 0$, we have the following identities
\begin{align*}
&\beta_{n,\lambda}(x+1-\lambda) =n(x-\lambda)_{n-1,\lambda}+\beta_{n,\lambda}(x-\lambda), \\
& \mathcal{E}_{n,\lambda}(x+1-\lambda) = 2(x-\lambda)_{n,\lambda}- \mathcal{E}_{n,\lambda}(x-\lambda), \\
&\mathcal{G}_{n,\lambda}(x+1-\lambda) = 2n(x-\lambda)_{n-1,\lambda} - \mathcal{G}_{n,\lambda}(x-\lambda).
\end{align*}

The study of degenerate versions of special polynomials and numbers originated with Carlitz's pioneering work on degenerate Bernoulli and Euler numbers [2,3]. This field has since expanded to include transcendental functions, such as gamma functions [16]. Furthermore, the introduction of $\lambda$-umbral calculus [13] has established a framework more robust than classical umbral calculus for analyzing degenerate Sheffer polynomials. These investigations employ a diverse array of methodologies, including generating functions, combinatorial methods, $p$-adic analysis, operator theory, differential equations, and probability theory (see [10-12,14-17,23] and references therein). For general background and references, one may consult [6,25].
\vspace{0.1in} \par

\section{Degenerate Bernoulli, Euler, and Genocchi Numbers}

From \eqref{3}, we note that
\begin{equation}\label{16}
\beta_{n,\lambda}(1) - \beta_{n,\lambda}=
\begin{cases}
1, & \text{if } n = 1, \\[4pt]
0, & \text{if } n \neq 1.
\end{cases}
\end{equation}
Then, by \eqref{4} and \eqref{16}, we get
\begin{equation}\label{17}
\sum_{k=0}^{n}\binom{n}{k}(1)_{n-k,\lambda} \beta_{k,\lambda}  - \beta_{n,\lambda}
= \delta_{1,n},
\end{equation}
where $\delta_{n,k}$ is the Kronecker's symbel.
From \eqref{17}, we have
\begin{align}
&\beta_{0,\lambda}=1,\quad
\beta_{1,\lambda}=-\frac{1}{2}+\frac{1}{2}\lambda,\quad
\beta_{2,\lambda}=\frac{1}{6}-\frac{1}{6}\lambda^{2},\quad
\beta_{3,\lambda}=-\frac{1}{4}+\frac{1}{4}\lambda^{3},\quad \label{18} \\
&\beta_{4,\lambda}=-\frac{1}{30}+\frac{2}{3}\lambda^{2}-\frac{19}{30}\lambda^{4},\quad
\beta_{5,\lambda}=\frac{1}{4}\lambda-\frac{5}{2}\lambda^{3}+\frac{9}{4}\lambda^{5}, \dots. \nonumber
\end{align}
From \eqref{5}, we see that
\begin{equation}\label{19}
\mathcal{E}_{n,\lambda}(1)+\mathcal{E}_{n,\lambda}=2\delta_{0,n}.
\end{equation}
Then, from \eqref{6} and \eqref{19}, we obtain
\begin{equation}\label{20}
\sum_{k=0}^{n}\binom{n}{k}(1)_{n-k,\lambda} \mathcal{E}_{k,\lambda}+\mathcal{E}_{n,\lambda}=2\delta_{0,n}.
\end{equation}
Thus, by \eqref{20}, we get
\begin{align}
&\mathcal{E}_{0,\lambda}=1,\quad
\mathcal{E}_{1,\lambda}=-\frac{1}{2},\quad
\mathcal{E}_{2,\lambda}=\frac{\lambda}{2},\quad
\mathcal{E}_{3,\lambda}=\frac{1}{4} - \frac{3}{4}\lambda^2, \quad
\mathcal{E}_{4,\lambda}=-\frac{3}{2}\lambda + \frac{3}{2}\lambda^3, \label{21} \\
&\mathcal{E}_{5,\lambda}=-\frac{1}{2} + \frac{15}{2}\lambda^2 - \frac{15}{4}\lambda^4, \quad
\mathcal{E}_{6,\lambda}=\frac{15}{2}\lambda - 45\lambda^3 + \frac{45}{4}\lambda^5, \dots. \nonumber
\end{align}

From \eqref{7}, we get
\begin{equation}\label{22}
\mathcal{G}_{n,\lambda}(1)+ \mathcal{G}_{n,\lambda}=2\delta_{1,n}, \quad \mathcal{G}_{0,\lambda}=0.
\end{equation}
Then, from \eqref{8} and \eqref{22}, we have
\begin{equation}\label{23}
\sum_{k=0}^{n}\binom{n}{k}(1)_{n-k,\lambda}\mathcal{G}_{k,\lambda}+\mathcal{G}_{n,\lambda}=2\delta_{1,n}.
\end{equation}
From (23), we note that
\begin{align}
&\mathcal{G}_{1,\lambda}=1,\quad
\mathcal{G}_{2,\lambda}=-1,\quad
\mathcal{G}_{3,\lambda}=\frac{3}{2}\lambda, \quad
\mathcal{G}_{4,\lambda}=1 - 3\lambda^2, \quad
\mathcal{G}_{5,\lambda}=-\frac{15}{2}\lambda + \frac{15}{2}\lambda^3,  \label{24} \\
&\mathcal{G}_{6,\lambda}=-3 + 45\lambda^2 - \frac{45}{2}\lambda^4, \quad
\mathcal{G}_{7,\lambda}=\frac{105}{2}\lambda - 315\lambda^3 + \frac{315}{4}\lambda^5,\dots. \nonumber
\end{align}

\section{Degenerate Euler--Seidel method for degenerate Bernoulli,
Euler, and Genocchi polynomials}

For a given sequence $\big(a_{n}(x|\lambda)\big)_{n \ge 0}$, we consider the degenerate Euler-Seidel matrix associated with this sequence which is determined recursively by (see \eqref{10})
\begin{equation}\label{25}
\begin{aligned}
&a_{0,n}(x\mid\lambda)=a_n(x\mid\lambda), \quad (n \ge 0), \\
&a_{k,n}(x\mid\lambda) = \bigl(1-(k-n)\lambda\bigr)\, a_{k-1,n}(x\mid\lambda)
+ a_{k-1,n+1}(x\mid\lambda), \quad (n \ge 0,\, k \ge 1).
\end{aligned}
\end{equation}
Using \eqref{25} and by induction, we get
\begin{equation}\label{26}
a_{n,0}{(x\mid\lambda)}=\sum_{k=0}^{n}\binom{n}{k}(1-\lambda)_{\,n-k,\lambda}\,a_{0,k}{(x\mid\lambda)},
\end{equation}
and
\begin{equation}\label{27}
a_{0,n}{(x\mid\lambda)}=\sum_{k=0}^{n}\binom{n}{k}(-1)^{\,n-k}\langle 1-\lambda\rangle_{\,n-k,\lambda}\,
a_{k,0}{(x\mid\lambda)},
\end{equation}
where the degenerate rising factorial sequence is given by
$$
\langle x\rangle_{0,\lambda}=1, \ \ \langle x\rangle_{n,\lambda}= x(x+\lambda)(x+2\lambda)\cdots(x+(n-1)\lambda),
\quad (n\ge1).$$
Therefore, by \eqref{26} and
\eqref{27}, we obtain the following theorem.

\begin{theorem}\label{thm:1}
For $n \ge 0$, we have
\begin{equation*}
a_{n,0}{(x\mid\lambda)}=\sum_{k=0}^{n}\binom{n}{k}(1-\lambda)_{\,n-k,\lambda}\,a_{0,k}{(x\mid\lambda)},
\end{equation*}
and
\begin{equation*}
a_{0,n}{(x\mid\lambda)}=\sum_{k=0}^{n}\binom{n}{k}(-1)^{\,n-k}\langle 1-\lambda\rangle_{\,n-k,\lambda}\,
a_{k,0}{(x\mid\lambda)}.
\end{equation*}
\end{theorem}

Let $A_\lambda(x,t)$ be the generating function of the sequence
$\big(a_{0,n}(x\mid\lambda)\big)_{n \ge 0}$,\, given by
\begin{equation}\label{28}
A_\lambda(x,t)=\sum_{n=0}^{\infty} a_{0,n}{(x\mid\lambda)}\frac{t^n}{n!},
\end{equation}
and let $\overline{A}_\lambda(x,t)$ be that of the sequence
$\big(a_{n,0}(x\mid\lambda)\big)_{n \ge 0}$,\, given by
\begin{equation}\label{29}
\overline{A}_\lambda(x,t)=\sum_{n=0}^{\infty} a_{n,0}{(x\mid\lambda)}\frac{t^n}{n!}.
\end{equation}
Then, by Theorem 3.1, we have
\begin{equation}
\begin{aligned}\label{30}
\overline{A}_\lambda(x,t)&=\sum_{n=0}^{\infty} a_{n,0}{(x\mid\lambda)}\frac{t^n}{n!}  \\
&=\sum_{n=0}^{\infty}\sum_{k=0}^{n}\binom{n}{k}(1-\lambda)_{n-k,\lambda}\, a_{0,k}{(x\mid\lambda)}\frac{t^n}{n!} \\
&=\sum_{k=0}^{\infty}\frac{1}{k!}  a_{0,k}{(x\mid\lambda)}
\sum_{n=k}^{\infty}\frac{(1-\lambda)_{n-k,\lambda}}{(n-k)!}
{t^n} \\
&=\sum_{k=0}^{\infty}a_{0,k}{(x\mid\lambda)}\frac{t^k}{k!}\sum_{n=0}^{\infty}(1-\lambda)_{n,\lambda}
\frac{t^n}{n!} \\
&=A_\lambda(x,t)\,e_\lambda^{\,1-\lambda}(t).
\end{aligned}
\end{equation}
Therefore, by \eqref{30}, we obtain the following theorem.

\begin{theorem}\label{thm:2}
Let
\begin{equation*}
A_\lambda(x,t)=\sum_{n=0}^{\infty} a_{0,n}{(x\mid\lambda)}\frac{t^n}{n!}.
\end{equation*}
Then we have
\begin{equation}\label{31}
\overline{A}_\lambda(x,t)=\sum_{n=0}^{\infty} a_{n,0}{(x\mid\lambda)}\frac{t^n}{n!}
=e_\lambda^{\,1-\lambda}(t)\,A_\lambda(x,t).
\end{equation}
\end{theorem}
Let $a_{0,n}{(x\mid\lambda)}=\beta_{n,\lambda}(x), \ (n\ge0)$. Then we have
\begin{equation}\label{32}
A_\lambda(x,t)=\sum_{n=0}^{\infty} a_{0,n}(x\mid\lambda)\frac{t^n}{n!}
=\sum_{n=0}^{\infty} \beta_{n,\lambda}(x)\frac{t^n}{n!} =\frac{t}{e_\lambda(t)-1}\,e_\lambda^{\,x}(t).
\end{equation}
From \eqref{31} and \eqref{32}, we note that
\begin{equation}\label{33}
\begin{aligned}
\sum_{n=0}^{\infty} a_{n,0}{(x\mid\lambda)} \frac{t^n}{n!}
&=e_\lambda^{\,1-\lambda}(t)\,A_\lambda(x,t) =\frac{t}{e_\lambda(t)-1}\,e_\lambda^{x+1-\lambda}(t)\\
&=\sum_{n=0}^{\infty} \beta_{n,\lambda}(x+1-\lambda)\frac{t^n}{n!}.
\end{aligned}
\end{equation}
Thus, by comparing the coefficients on both sides of \eqref{33}, we get
\begin{equation}\label{34}
a_{n,0}{(x\mid\lambda)}=\beta_{n,\lambda}(x+1-\lambda),\quad (n\ge0).
\end{equation}
On the other hand, by \eqref{3}, \eqref{31} and \eqref{32}, we get
\begin{equation}\label{35}
\begin{aligned}
\sum_{n=0}^{\infty} a_{n,0}{(x\mid\lambda)}\frac{t^n}{n!}
&=e_\lambda^{\,1-\lambda}(t)\,A_\lambda(x,t)
=e_\lambda(t)\,\frac{t}{e_\lambda(t)-1}\,e_\lambda^{x-\lambda}(t) \\
&=(e_\lambda(t)-1+1)\frac{t}{e_\lambda(t)-1}\,e_\lambda^{\,x-\lambda}(t)\\
& =t e_\lambda^{x-\lambda}(t) + \frac{t}{e_\lambda(t)-1} e_\lambda^{\,x-\lambda}(t) \\
&=\sum_{n=0}^{\infty}\left(n(x-\lambda)_{n-1,\lambda}+\beta_{n,\lambda}(x-\lambda)\right)\frac{t^n}{n!}.
\end{aligned}
\end{equation}
Comparing the coefficients on both sides of
\eqref{35}, we have
\begin{equation}\label{36}
a_{n,0}{(x\mid\lambda)}=n(x-\lambda)_{n-1,\lambda}+\beta_{n,\lambda}(x-\lambda),\quad (n\ge0).
\end{equation}
Therefore, by \eqref{34} and \eqref{36}, we obtain the following theorem.

\begin{theorem}\label{thm:3}
For $n \ge 0$, we have
\begin{equation*}
\beta_{n,\lambda}(x+1-\lambda) =n(x-\lambda)_{n-1,\lambda}+\beta_{n,\lambda}(x-\lambda).
\end{equation*}
\end{theorem}

The degenerate Euler--Seidel matrix associated with $\left(a_{0,n}(x\mid\lambda)\right)_{n\ge 0}=(\beta_{n,\lambda}(x))_{n\ge0}$ is given by (see \eqref{3},\eqref{23},
\eqref{36})

\begin{equation*}
\left(\begin{matrix}
1 &  x - \frac{1-\lambda}{2} & \cdots \\
x + \frac{1-\lambda}{2} & x^2 - \frac{\lambda^2 - 3\lambda + 2}{6} & \cdots \\
x^2 + (1-2\lambda)x + \frac{1-6\lambda+5\lambda^2}{6}
&  x^3 + \frac{1-3\lambda}{2}x^2
  + \frac{6\lambda^2 - 5\lambda - 1}{2}x
  + \frac{4\lambda^3 - 3\lambda^2 - 1}{6} & \cdots \\
\vdots & \vdots & \vdots
\end{matrix}\right).
\end{equation*}

Let $a_{0,n}(x) = \mathcal{E}_{n,\lambda}(x), (n \ge 0)$. Then we have
\begin{equation} \label{37}
A_\lambda(x,t)=\sum_{n=0}^{\infty} a_{0,n}(x\mid\lambda)\,\frac{t^n}{n!}
=\sum_{n=0}^{\infty} \mathcal{E}_{n,\lambda}(x)\,\frac{t^n}{n!}
=\frac{2}{e_\lambda(t) + 1} \, e_\lambda^{x}(t).
\end{equation}
From \eqref{31} and \eqref{37}, we note that
\begin{equation}
\begin{split} \label{38}
\sum_{n=0}^{\infty} a_{n,0}(x\mid\lambda)\,\frac{t^n}{n!}&= \overline{A}_{\lambda}(x,t)
=e_\lambda^{1-\lambda}(t) A_\lambda(x,t)=e_\lambda^{1-\lambda}(t) \frac{2}{e_\lambda(t)+1}  {e_\lambda^x (t)}  \\
& = \frac{2}{e_\lambda(t)+1}  {e_\lambda^{x+1-\lambda} (t)} = \sum_{n=0}^{\infty}\mathcal{ E}_{n,\lambda}(x+1-\lambda)\,\frac{t^n}{n!}.
\end{split}
\end{equation}
Thus, by \eqref{38}, we get
\begin{equation}\label{39}
a_{n,0}(x \mid \lambda)=\mathcal{E}_{n,\lambda}(x+1-\lambda),\ \ (n \ge 0).
\end{equation}
On the other hand, by \eqref{4} and \eqref{31}, we get
 \begin{equation}\label{40}
 \begin{split}
&\sum_{n=0}^{\infty} a_{n,0}(x \mid \lambda)\,\frac{t^n}{n!}=
\overline{A}_\lambda(x,t)=e_\lambda^{1-\lambda}(t) A_\lambda(x,t)   \\
&=\left(e_\lambda(t)+1-1\right)\frac{2}{e_\lambda(t)+1}e_\lambda^{x-\lambda}(t)=2e_\lambda^{x-\lambda}(t)
-  \frac{2}{e_\lambda(t)+1}e_\lambda^{x-\lambda}(t)  \\
&= \sum_{n=0}^{\infty} \left(2(x-\lambda)_{n,\lambda} - \mathcal{E}_{n,\lambda}(x-\lambda)\right)\frac{t^n}{n!}.
 \end{split}
 \end{equation}
Comparing the coefficients on both sides of \eqref{40}, we have
 \begin{equation}\label{41}
a_{n,0}(x \mid \lambda)=2(x-\lambda)_{n,\lambda}-\mathcal{E}_{n,\lambda}(x-\lambda),\quad (n \ge 0).
\end{equation}
Therefore, by \eqref{39} and \eqref{41}, we obtain the following theorem.

\begin{theorem}\label{thm:4} For  $ n \ge 0$,  we have
$$ \mathcal{E}_{n,\lambda}(x+1-\lambda) = 2(x-\lambda)_{n,\lambda}- \mathcal{E}_{n,\lambda}(x-\lambda). $$
\end{theorem}

By \eqref{5}, \eqref{23}, and \eqref{41}, we see that the degenerate Euler--Seidel matrix associated with
$\big(a_{0,n}\big)_{n\ge 0} = \big(\mathcal{E}_{n,\lambda}(x)\big)_{n\ge 0}$
 is given by

\begin{equation}
 \left(\begin{matrix}
 1 & x-\dfrac12 &  \cdots\\
x+\dfrac12-\lambda  & x^2-\lambda x+\dfrac{\lambda-1}{2} & \cdots\\
x^2+(1-3\lambda)x+\dfrac{4\lambda^2-3\lambda}{2}  & x^3+\Bigl(\dfrac12+\lambda\Bigr)x^2+(2\lambda^2-1)x-\dfrac14(1-6\lambda+2\lambda^2) &  \cdots\\
 \vdots & \vdots & \vdots
 \end{matrix}\right).\label{42}
 \end{equation}

Finally, we let $a_{0,n}(x \mid \lambda) = \mathcal{G}_{n,\lambda}(x)$, ($n \ge 0$).
Then, by \eqref{7} and \eqref{28}, we get
\begin{equation}\label{43}
A_\lambda(x,t)= \sum_{n=0}^{\infty} a_{0,n}(x \mid \lambda)\,\frac{t^n}{n!}
=\sum_{n=0}^{\infty} \mathcal{G}_{n,\lambda}(x)\,\frac{t^n}{n!}=
\frac{2t}{e_\lambda{(t)}+1}\, e_\lambda^{x}(t).
\end{equation}
From \eqref{31} and \eqref{43}, we note that
\begin{equation}\label{44}
\begin{aligned}
\sum_{n=0}^{\infty} a_{n,0}{(x \mid \lambda)}\frac{t^n}{n!}
&= \overline{ A}_\lambda(x,t) = e_\lambda^{\,1-\lambda}(t)\,A_\lambda(x,t) =e_\lambda^{\,1-\lambda}(t) \frac{2t}{e_\lambda(t)+1}\,e_\lambda^{\,x}(t) \\
&= \frac{2t}{e_\lambda(t)+1}\,e_\lambda^{\,x+1-\lambda}(t)
=\sum_{n=0}^{\infty} \mathcal{G}_{n,\lambda}(x+1-\lambda)\frac{t^n}{n!}.
\end{aligned}
\end{equation}
Thus, by \eqref{44}, we get
\begin{equation}\label{45}
\mathcal{G}_{n,\lambda}(x+1-\lambda)= a_{n,0}{(x \mid \lambda)},\quad (n\ge0).
\end{equation}
On the other hand, by \eqref{7}, \eqref{31}, and \eqref{43}, we get
\begin{equation}\label{46}
\begin{aligned}
&\sum_{n=0}^\infty a_{n,0}(x \mid \lambda) \frac{t^n}{n!}
=\overline{A}_\lambda(x,t) = e_\lambda^{\,1-\lambda}(t)\,A_\lambda(x,t) \\
&=\frac{2t}{e_\lambda(t)+1}\,e_\lambda^{\,x-\lambda}(t)\left(e_\lambda(t)+1-1\right)
={2t}e_\lambda^{x-\lambda}(t)-\frac{2t}{e_\lambda(t)+1}\,e_\lambda^{\,x-\lambda}(t) \\
&= \sum_{n=0}^{\infty} \Bigl(2n(x-\lambda)_{n-1,\lambda} - \mathcal{G}_{n,\lambda}(x-\lambda)
\Bigr) \frac{t^n}{n!}.
\end{aligned}
\end{equation}
By comparing the coefficients on both sides of \eqref{46}, we get
\begin{equation}\label{47}
a_{n,0}{(x \mid \lambda)}= 2n(x-\lambda)_{n-1,\lambda} - \mathcal{G}_{n,\lambda}(x-\lambda), \qquad (n\ge0).
\end{equation}
Therefore, by \eqref{45} and \eqref{47}, we obtain the following theorem.

\begin{theorem}\label{thm:5}
For $n \ge 0$, we have
\begin{equation*}
\mathcal{G}_{n,\lambda}(x+1-\lambda) = 2n(x-\lambda)_{n-1,\lambda} - \mathcal{G}_{n,\lambda}(x-\lambda).
\end{equation*}
\end{theorem}

From \eqref{8}, \eqref{25} and \eqref{47}, we note that the degenerate Euler--Seidel
matrix associated with the sequence
$\left(a_{n,0}(x \mid \lambda)\right)_{n\ge 0}=\left(\mathcal{G}_{n,\lambda}(x)\right)_{n\ge0}$ is given by

 \begin{equation}
\left(\begin{matrix}
 0 & 1 & 2x-1 & \cdots\\
1  &2x  & \cdots  &  \cdots\\
 2x+1-2\lambda &  3x^{2}+(1-3\lambda)x-1+\dfrac{\lambda}{2} & \cdots & \cdots\\
3x^{2}+3(1-3\lambda)x+6\lambda^{2}-\dfrac{9\lambda}{2}& \cdots & \cdots & \cdots\\
 \vdots & \vdots & \vdots & \vdots
 \end{matrix}\right).\label{48}
 \end{equation}
\medskip

\begin{remark}
Note that
\[
\lim_{\lambda \to 0} \beta_{n,\lambda}(x) = B_n(x), \qquad
\lim_{\lambda \to 0} \mathcal{E}_{n,\lambda}(x) = {E}_n(x), \qquad
\lim_{\lambda \to 0} \mathcal{G}_{n,\lambda}(x) = G_n(x),
\]
where $B_n(x)$, $E_n(x)$, and $G_n(x)$ are the ordinary Bernoulli, Euler,
and Genocchi polynomials given by (see [4,7,21,26,27])
$$ \frac{t}{e^t-1} e^{xt}=\sum_{n=0}^{\infty} B_n(x)\frac{t^n}{n!},$$
$$ \frac{2}{e^t+1} e^{xt}=\sum_{n=0}^{\infty} {E}_n(x)\frac{t^n}{n!},$$
and
$$\frac{2t}{e^t+1} e^{xt}= \sum_{n=0}^{\infty} G_n(x)\frac{t^n}{n!}.$$
\end{remark}

\bigskip

\section{Conclusion}
In this study, we successfully generalized the Euler-Seidel method to its degenerate form. We demonstrated that the transition from an initial sequence to a final sequence in the degenerate Euler-Seidel matrix can be expressed through $\lambda$-generalized binomial identities involving degenerate falling and rising factorials. Furthermore, we showed that the relationship between their respective generating functions is governed by the degenerate exponential function $e_{\lambda}^{1-\lambda}(t)$. By implementing this method, we obtained explicit identities for degenerate Bernoulli, Euler, and Genocchi polynomials.
These results not only recover classical identities as $\lambda \rightarrow 0$ but also offer a robust analytical tool for investigating broader classes of degenerate special sequences in combinatorial analysis and number theory.


\end{document}